\documentclass[a4papersize]{article}
\usepackage{amsmath,amsthm, amsxtra,amssymb,latexsym, amscd}
\usepackage[mathscr]{eucal}
\usepackage{color}
\newtheorem{theorem}{Theorem}[section]

\newtheorem{lemma}[theorem]{Lemma}
\newtheorem{proposition}[theorem]{Proposition}
\newtheorem{example}[theorem]{Example}

\newtheorem{remark}[theorem]{Remark}

\DeclareMathOperator{\Hom}{Hom}
\DeclareMathOperator{\Ann}{Ann}

\DeclareMathOperator{\Ass}{Ass}
\DeclareMathOperator{\Assh}{Assh}
\DeclareMathOperator{\Supp}{Supp}

\DeclareMathOperator{\Var}{Var}

\DeclareMathOperator{\nCM}{nCM}

\DeclareMathOperator{\Rad}{Rad}
\DeclareMathOperator{\Att}{Att}

\DeclareMathOperator{\Spec}{Spec}
\DeclareMathOperator{\p}{\frak p}

\DeclareMathOperator{\Ext}{Ext}
\begin{document}
\large
\centerline{\Large {\bf ANNIHILATOR OF LOCAL COHOMOLOGY MODULES}}
\medskip
\centerline{\Large {\bf UNDER LOCALIZATION AND COMPLETION}}

\vskip 0.7cm
\vskip 0.4cm\centerline { NGUYEN THI ANH HANG}
\centerline {Thai Nguyen University of Education}
\centerline{Thai Nguyen, Vietnam}
\centerline {e-mail: hangnthianh@gmail.com}
\vskip 0.4cm

\centerline {LE THANH NHAN}
\centerline {Institute of Mathematics}
\centerline{Vietnam Academy of Science and Technology, Hanoi, Vietnam}
\centerline {E-mail: nhanlt2014@gmail.com}

\vskip 0.4cm
\smallskip
\vskip 1cm

\noindent{\bf Abstract} {\footnote{ {\it{Key words and phrases: }} Annihilator of local cohomology module, catenarity, unmixedness, Cohen-Macaulayness. \hfill\break
  {\it{2020 Subject  Classification: }} 13D45, 13H10, 13E05.\hfill\break {This research is funded by Vietnam Ministry of Education and Training  under grant number B2025-TNA-03.}}. {Let $(R, \frak m)$ be a Noetherian local ring. This paper deals with the annihilator of Artinian local cohomology modules $H^i_{\frak m}(M)$ in the relation with the structure of the base ring $R$, for non negative integers $i$ and finitely generated $R$-modules $M$. Firstly, the catenarity and the unmixedness of local rings are characterized via the compatibility of annihilator of top local cohomology modules  under localization and completion, respectively. Secondly, some necessary and sufficient conditions for a local ring being a quotient of a Cohen-Macaulay local ring are given in term of the annihilator of all local cohomology modules under localization and completion.}

\section{Introduction}

\ \ \ \  \ \ Throughout this paper,  let $(R,\frak m )$ be a Noetherian local ring. 

It is well known that the annihilator of finitely generated $R$-modules is compatible under localization, i.e. $\Ann_{R_{\frak p}}M_{\frak p}=\big(\Ann_RM)R_{\frak p}$  for every  finitely generated $R$-module $M$ and every $\frak p\in\Spec(R)$. Note that $H^i_{\frak m}(M)$ is an Artinian $R$-module  for any finitely generated $R$-module $M$ and any integer $i\geq 0$. Therefore, it is natural to ask whether  the annihilator of $H^i_{\frak m}(M)$ is compatible under localization, i.e. 
$$\Ann_{R_{\frak p}}\big(H^{i-\dim (R/\frak p)}_{\frak p R_{\frak p}}(M_{\frak p})\big)=\big(\Ann_{R}H^{i}_{\frak m}(M)\big)R_{\frak p}$$ for every finitely generated $R$-module $M$,  every integer $i \geq 0$ and every $\frak p \in \Spec(R)$. In general, the answer is negative, see Examples \ref{E:1}, \ref{E:2}. 

The annihilator of finitely generated modules is also compatible under $\frak m$-adic completion, i.e. $\Ann_{\widehat R}\widehat M=(\Ann_RM)\widehat R$ for every finitely generated $R$-module $M$. Note that $H^i_{\frak m}(M)$ is an Artinian module over both $R$ and $\widehat R$. Moreover, $H^i_{\frak m}(M)\otimes_R\widehat R\cong H^i_{\frak m}(M)$. So, we ask whether the annihilator of $H^i_{\frak m}(M)$ is compatible under completion, i.e. 
$$\Ann_{\widehat R}H^i_{\frak m}(M)=(\Ann_R H^i_{\frak m}(M))\widehat R$$ for every finitely generated $R$-module $M$,  every integer $i \geq 0$. In general, the answer is negative, see Examples \ref{E:1}, \ref{E:2}.  

The aim of this paper is to study the structure of the base ring $(R,\frak m)$ satisfying the condition that the annihilator of local cohomology modules with support in $\frak m$ is compatible under localization and completion.

The first main result of this paper gives a characterization of catenary local rings via the compatibility of annihilator of the top local cohomology modules under localization.
  
\begin{theorem} \label{T:1a} The following statements are equivalent:

(i) $R$ is catenary.

(ii) $\Ann_{R_{\frak p}}\big(H^{\dim_R(M)-\dim (R/\frak p)}_{\frak p R_{\frak p}}(M_{\frak p})\big)=\big(\Ann_{R}H^{\dim_R(M)}_{\frak m}(M)\big)R_{\frak p}$ for every finitely generated $R$-module $M$ of dimension $d_M$ and every $\frak p \in \Spec(R)$.
\end{theorem}

The second main result of this paper clarifies the structure of local rings via the compatibility of annihilator of the top local cohomology modules under completion. Following M. Nagata \cite{Na}, $M$ is said to be {\it unmixed}  if $\dim(\widehat R/\frak P)=\dim_R(M)$ for all $\frak P\in\Ass_{\widehat R} (\widehat M)$.

\begin{theorem} \label{T:2a} The following statements are equivalent:

(i) $R/\frak p$ is unmixed for any $\frak p \in \Spec(R)$.
	
(ii) $\Ann_{\widehat R} H_{\frak m}^{\dim_R(M)}(M)=(\Ann_R H_{\frak m}^{\dim_R(M)}(M)) \widehat R$ for any finitely generated $R$-module $M$.
\end{theorem}

Theorems \ref{T:1a}, \ref{T:2a} are concerned with the annihilator of the top local cohomology modules. Now we consider the annihilator of all local cohomology modules $H^i_{\frak m}(M)$ for finitely generated $R$-modules $M$ and integers $i\leq \dim_R(M)$.  If $R$ is a quotient of a Gorenstein local ring, we can use Local Duality Theorem to  show that the annihilator of all local cohomology modules $H^i_{\frak m}(M)$ is compatible under completion and localization, see Proposition \ref{P:1}. Note that there exists a local ring $(R,\frak m)$ which is not a quotient of a Gorenstein local ring, but the annihilator of all local cohomology modules is still compatible under completion and localization, see Remark \ref{R:3}. 

The third main result of this paper gives characterizations for the base ring being a quotient of a Cohen-Macaulay local ring via the compatibility of annihilator of all local cohomology modules under localization and completion. 

\begin{theorem}\label{T:3a} The following statements are equivalent:

(i) $R$ is a quotient of a Cohen-Macaulay local ring.

(ii)  $\Rad \Big(\Ann_{R_{\frak p}} H^{i-\dim (R/\frak p)}_{\frak p R_{\frak p}}(M_{\frak p})\Big) = \Rad\Big(\big(\Ann_R H^i_{\frak m}(M)\big)R_{\frak p}\Big)$ for every finitely generated $R$-module $M$, every $\frak p \in \Supp_R(M)$ and every $i \geq 0$.
	
(iii) $ \Rad \big(\Ann_{\widehat{R}} H^i_{\frak m}(M)\big) = \Rad\Big(\big(\Ann_{R} H^i_{\frak m}(M)\big)\widehat{R}\Big)$ for every finitely generated $R$-module $M$ and every integer $i\geq 0$. 
\end{theorem}

The proofs of Theorems \ref{T:1a}, \ref{T:2a} will be presented in the next section.  In the last section, we prove Theorem \ref{T:3a} and give some examples to clarify the results of this paper.

\section{Proof of Theorems \ref{T:1a} and \ref{T:2a}}

Throughout this paper, we use the following notations: For each finitely generated $R$-module $M$ of dimension $d$,  we set  $$\Assh_R(M)=\{\frak p\in\Ass_R(M)\mid \dim(R/\frak p)=d\}.$$ 
For a submodule $N$ of $M$ and a reduced primary decomposition $N=\bigcap\limits_{\frak p\in\Ass_R(M/N)}N(\frak p)$ of $N$, we set $U_M(N):= \bigcap\limits_{\frak p\in\Assh_R(M/N)}N(\frak p).$ Note that $U_M(N)$  does not depend on the choice of reduced primary decomposition of $N$. Moreover, $N\subseteq U_M(N)$ and $U_M(N)/N$ is the largest submodule of $M/N$ of dimension less than $\dim_R(M/N).$

The following result can be seen in \cite[Corollary 1.3]{ASN}, \cite[Theorem 2.6]{BAG}.

\begin{lemma} \label{L:1a}  $\Ann_R H^d_{\frak m}(M)=\Ann_R (M/U_M(0))$.
\end{lemma}

We need to recall some facts of attached primes for Artinian modules. Following I. G. Macdonald \cite{Mac}, it is well known that every Artinian $R$-module $A$ has a minimal secondary representation $A=A_1+\ldots +A_n$, where each $A_i$ is $\frak p_i$-secondary,  $A_i$ is not redundant and $\frak p_i\neq \frak p_j$ for all $i\neq j.$ The set $\{\frak p_1, \ldots , \frak p_n\}$ does not depend on the choice of minimal secondary representation of $A$. This set is denoted by $\Att_RA$ and it is called the  set of {\it attached primes} of $A$. For an ideal $I$ of $R$, denote by $\Var(I)$ the set of all prime ideals of $R$ containing $I$. The following facts can be seen in \cite{BS}, \cite{Mac}.

\begin{lemma}\label{L:1} Let $A$ be an Artinian $R$-module.  Let $M$ be a finitely generated $R$-module of dimension $d$. Then 

(i) $A\neq 0$ if and only if $\Att_RA\neq \emptyset.$ Moreover, $\ell_R(A)<\infty$ if and only if $\Att_RA \subseteq \{\frak m\}.$

(ii) $\min \Att_RA=\min \Var (\Ann_RA).$ 

(iii) $A$ has a natural structure as an Artinian $\widehat R$-module and $$\Att_RA=\{\frak P\cap R\mid \frak P\in \Att_{\widehat R}A\}.$$

(iv) $\Att_RH^d_{\frak m}(M)=\Assh_R(M).$
\end{lemma}
  \begin{proof}[{\bf Proof of Theorem \ref{T:1a}}] (i) $\Rightarrow$ (ii). Let $M$ be a finitely generated $R$-module of dimension $d$ and let $\frak p\in \Spec(R)$. Suppose that $\frak p\notin \Supp_R(M/U_M(0))$. Then $$\dim_{R_{\frak p}}(M_{\frak p})=\dim_{R_{\frak p}}(U_M(0)_{\frak p})\leq \dim_R(U_M(0))-\dim (R/\frak p)<d-\dim (R/\frak p).$$ Hence $H^{d-\dim (R/\frak p)}_{\frak p R_{\frak p}}(M_{\frak p})=0$ by \cite[Theorem 6.1.4]{BS}. So, $\Ann_{R_{\frak p}}H^{d-\dim (R/\frak p)}_{\frak p R_{\frak p}}(M_{\frak p})=R_{\frak p}.$ Moreover, we have by Lemma \ref{L:1a} that 
$$(\Ann_{R}H^d_{\frak m}(M))R_{\frak p}=(\Ann_RM/U_M(0))R_{\frak p}=R_{\frak p}.$$
So, we can assume that $\frak p\in \Supp_R(M/U_M(0))$. Then $\frak p\supseteq \frak p'$ for some $\frak p'\in \Assh_R(M).$ Since $R$ is catenary, 
$$d=\dim (R/\frak p')=\dim(R/\frak p)+{\rm ht}(\frak p/\frak p')\leq \dim(R/\frak p)+\dim(M_{\frak p})\leq d.$$
Hence $\dim_{R_{\frak p}} (M_{\frak p})=d-\dim (R/\frak p).$ 	Let
$0=\bigcap\limits_{\frak q \in \Ass_R(M)} N(\frak q)$ be a reduced primary decomposition of $0$ in $M$.	We can write $\Assh_R(M)=\{\frak q_1, \ldots, \frak q_r\}$, where $0<s\leq r$ is an integer such that $\frak q_i \subseteq \frak p$ for $i=1, \ldots, s$ and $\frak q_i \not\subseteq \frak p$ for $i=s+1, \ldots, r.$ Since $U_M(0)=\bigcap\limits_{i=1}^r N(\frak q_i)$, we have $U_M(0)_{\frak p}=\bigcap\limits_{i=1}^s N(\frak q_i)_{\frak p}.$ Now let $\frak q \in \Ass_R(M)$ such that $\frak q \subseteq \frak p.$ 	As $R$ is catenary,   $$\dim (R_{\frak p}/\frak q R_{\frak p})= {\rm ht}(\frak p/\frak q)=\dim (R/\frak q)-\dim (R/\frak p).$$ Hence $\dim (R_{\frak p}/\frak q R_{\frak p})=\dim_{R_{\frak p}} (M_{\frak p})$ if and only if $\dim (R/\frak q)=d.$ Therefore, 
  	$$\Assh_{R_{\frak p}}(M_{\frak p})=\{\frak q_1 R_{\frak p}, \ldots, \frak q_s R_{\frak p}\}.$$
  	Note that $0_{M_{\p}}=\underset{\frak q \in \Ass_R(M), \frak q \subseteq \frak p}{\bigcap} N(\frak q)_{\frak p}$ is a reduced primary decomposition of $0$ in $M_{\frak p}$. So, $$U_{M_{\frak p}}(0)=\bigcap_{i=1}^s N(\frak q_i)_{\frak p}=U_M(0)_{\frak p}.$$
  	Therefore, we get by Lemma \ref{L:1a} that  	
  	\begin{align*}\Ann_{R_{\frak p}} H^{d-\dim (R/\frak p)}_{\frak p R_{\frak p}}(M_{\frak p}) &=\Ann_{R_{\frak p}} M_{\frak p}/U_{M_{\frak p}}(0) 
  		=\Ann_{R_{\frak p}}(M/U_M(0))_{\frak p}\\&=\big(\Ann_R(M/U_M(0)\big)R_{\frak p}=\big(\Ann_{R}H^{d}_{\frak m}(M)\big)R_{\frak p}.
  	\end{align*}
  	
  	(ii) $\Rightarrow$ (i). Let $\frak q\in \Spec(R)$. We use Ratliff's criterion \cite[Theorem 2.2]{Rat} to prove the catenarity of $R/\frak q$. Let $\frak p \in \Spec(R)$ such that $\frak q \subseteq \frak p$. Set $t=\dim (R/\frak q)$. Since  $\Att_R(H^{t}_{\frak m}(R/\frak q))=\{\frak q\}$,  we have $\frak q =\Rad(\Ann_{R}H^{t}_{\frak m}(R/\frak q))$ by Lemma \ref{L:1}(ii),(iv). Hence $\frak q =\Ann_{R}H^{t}_{\frak m}(R/\frak q)$. Therefore, we get by our assumption (ii) that
  	$$\Ann_{R_{\frak p}}\big(H^{t-\dim (R/\frak p)}_{\frak p R_{\frak p}}(R_{\frak p}/\frak q R_{\frak p})\big)=\frak q R_{\frak p}.$$ Hence $H^{t-\dim (R/\frak p)}_{\frak p R_{\frak p}}(R_{\frak p}/\frak q R_{\frak p})\neq 0$. It follows from \cite[Theorem 6.1.4]{BS} that $$\dim (R/\frak p)+\dim(R_{\frak p}/\frak q R_{\frak p})=\dim (R/\frak q).$$ Hence $R/\frak q$ is catenary for any $\frak q\in \Spec(R)$. Thus $R$ is catenary.
  \end{proof}

\begin{proof}[{\bf Proof of Theorem \ref{T:2a}}]
(i) $\Rightarrow$ (ii). Let $M$ be a finitely generated $R$-module of dimension $d.$ Note that $\widehat{U_M(0)}\subseteq U_{\widehat M}(0)$. We  claim that  $\widehat{U_M(0)}=U_{\widehat M}(0)$. Assume in contrary that $\widehat{U_M(0)}\neq U_{\widehat M}(0)$. Then there exists $\frak P\in \Ass_{\widehat R}(U_{\widehat M}(0)/\widehat{U_M(0)}).$ It follows that  $\frak P\in\Ass_{\widehat R}(\widehat M/\widehat{U_M(0)})$ and $\dim (\widehat R/\frak P)\leq \dim_{\widehat R}\big(U_{\widehat M}(0)/\widehat{U_M(0)}\big) <d$. Set $\frak p=\frak P\cap R$. Then we get by \cite[Theorem 23.2]{Mat} that  $\frak P\in\Ass (\widehat R/\frak p\widehat R)$ and $\frak p \in \Ass_R(M/U_M(0))$. Hence $\dim (R/\frak p) =d$. As $R/\frak p$ is unmixed by our assumption (i), $\dim (\widehat R/\frak P)=d$.  This gives a contradiction. So, the claim follows. Now, we get by Lemma \ref{L:1a} and by the claim that
\begin{align}(\Ann_R H_{\frak m}^d(M)) \widehat R =\big(\Ann_R M/U_M(0)\big)\widehat R&=\Ann_{\widehat R} \big(\widehat M/\widehat{U_M(0)}\big)\notag\\&=\Ann_{\widehat R}\big(\widehat M/U_{\widehat M}(0)\big) =\Ann_{\widehat R} H_{\frak m}^{d}(M).\notag\end{align}

(ii) $\Rightarrow$ (i). Let $\frak p\in \Spec(R)$. Set $t=\dim (R/\frak p)$. It follows from Lemma \ref{L:1}(ii), (iv) that 
$\Ann_R H_{\frak m}^t(R/\frak p)=\frak p.$ By the assumption (ii), we have by Lemma \ref{L:1a} that
$$\frak p \widehat R=\Ann_{\widehat R} H_{\frak m}^t(R/\frak p)=\Ann_{\widehat R} H_{\frak m\widehat R}^t(\widehat R/\frak p\widehat R)=U_{\widehat R}(\frak p \widehat R).$$
Hence $\Ass (\widehat R/\frak p \widehat R)=\Assh (\widehat R/\frak p \widehat R),$ i.e. $R/\frak p$ is unmixed.	
\end{proof}

\section{Proof of Theorem \ref{T:3a}}

Let $E(R/\frak m)$ be the injective hull of the residue field $R/\frak m$, let $D_R(-):=\Hom_R(-, E(R/\frak m))$ be the Matlis duality functor.

\begin{lemma} \label{L:1b} Let $\frak p\in\Spec(R)$. Let $A$ be an Artinian $R$-module, let $M$ be a finitely generated $R$-module. If $A=D_R(M)$, then $(\Ann_RA)\widehat R = \Ann_{\widehat R}A$ and
$\Ann_{R_{\frak p}}D_{R_{\frak p}}(M_{\frak p})=(\Ann_RA)_{\frak p}.$
\end{lemma}

\begin{proof}   Since $M$ is finitely generated,  $\Ann_{\widehat R}\widehat M=(\Ann_RM)\widehat R$  and $\Ann_{R_{\frak p}}M_{\frak p}=(\Ann_RM)_{\frak p}.$ Therefore, 
$$(\Ann_RA)\widehat R=(\Ann_RD_R(M))\widehat R=(\Ann_RM)\widehat R=\Ann_{\widehat R} \widehat M.$$
Note that $\widehat M\cong D_R(D_R(M))$ as $\widehat R$-modules. Since $D_R(M)$ is an Artinian $R$-module, it is an Artinian $\widehat R$-module and $D_R(D_R(M))=D_{\widehat R}(D_R(M)).$ So,  
$$\Ann_{\widehat R}\widehat M=\Ann_{\widehat R}D_{\widehat R}(D_R(M))=\Ann_{\widehat R}D_{\widehat R}(A)=\Ann_{\widehat R}A.$$
Moreover, 
$$
\Ann_{R_{\frak p}}D_{R_{\frak p}}(M_{\frak p})=\Ann_{R_{\frak p}}M_{\frak p}=(\Ann_RM)_{\frak p}\\
=(\Ann_RD_R(M))_{\frak p}=(\Ann_RA)_{\frak p}.$$
\end{proof}

In case where $R$ is a quotient of a Gorenstein local ring, by using Lemma \ref{L:1b} and Local Duality Theorem (see \cite[Theorem 11.2.6]{BS}), we can show that the annihilator of local cohomology modules is compatible under localization and completion.

  \begin{proposition}\label{P:1} Let $M$ be a finitely generated $R$-module. Assume that $R$ is a quotient of a Gorenstein local ring. Then
  	
  (i) $ \Ann_{\widehat R}H^i_{\frak m}(M) = (\Ann_RH^i_{\frak m}(M))\widehat R$ for all integers $i \geq 0$.
  
  (ii) $\Ann_{R_{\frak p}}H^{i-\dim (R/\frak p)}_{\frak p R_{\frak p}}(M_{\frak p})=(\Ann_{R}H^{i}_{\frak m}(M))_{\frak p}$ for all $\frak p \in \Spec(R)$ and all integers $i \geq 0$.
  \end{proposition}

  \begin{proof}  Let $(R',\frak m')$ be a Gorenstein local ring of dimension $n$ such that $R$ is a quotient of $R'$.    For each $i\in \mathbb N$, denote by
  	$K^i_{M}:=\Ext^{n'-i}_{R'}(M,R')$ the $i$-th defficiency module of $M$. For each integer  $i \ge 0$, note that $K^i_{M}$ is a finitely generated $R$-module and Local Duality Theorem  gives an isomorphism 
$H^{i}_{\frak m}(M)\cong D_R(K^i_M).$  Now, the statement (i) follows by Lemma \ref{L:1b}.

    	Let $\frak p \in \Spec(R)$. For each integer $i\geq 0$, we have an $R_{\frak p}$-isomorphism $(K^{i}_{M})_{\frak p}\cong K^{i-\dim (R/\frak p)}_{M_{\frak p}}$	and Local Duality Theorem gives an $R_{\frak p}$-isomorphism
  $$H^{i-\dim (R/\frak p)}_{\frak p R_{\frak p}}(M_{\frak p})\cong D_{R_{\frak p}}(K^{i-\dim (R/\frak p)}_{M_{\frak p}}).$$
  Therefore, the statement (ii) follows by Lemma \ref{L:1b}.
    \end{proof}

 \begin{remark} \label{R:3} {\rm If $\dim(R)=1$ then the annihilator of all local cohomology modules $H^i_{\frak m}(M)$ is compatible under localization and completion. Indeed, the case $i=0$ is obvious. As $\dim(R)=1$,  $R$ is catenary and  $R/\frak p$ is unmixed for all $\frak p\in\Spec(R)$. So, the case $i=1$ follows by Theorems \ref{T:1a}, \ref{T:2a}. Note that there is a Noetherian local ring of dimension $1$ which can not be expressed as  a quotient of a Gorenstein local ring, see \cite{FR}.}
\end{remark}

 From Proposition \ref{P:1} and Remark \ref{R:3}, it is natural to ask about the problem in a more general situation where $R$ is a quotient of a Cohen-Macaulay local ring. A partial answer to this question is given in Theorem \ref{T:3a}. Before proving this theorem, we give a criterion for the Cohen-Macaulayness of the formal fiber of $R$ at a prime ideal $\frak p$. 

\begin{proposition}\label{P:1c} Let $\frak p$ be a prime ideal of $R$.  If  $$\Rad \big({\mathrm{Ann}_{\widehat{R}} H^i_{\frak m}(R/\frak p)\big) = \Rad\big((\Ann_{R} H^i_{\frak m}(R/\frak p))\widehat{R}}\big)$$ for all integers $i\leq \dim(R/\frak p)$, then $R/\frak p$ is quasi unmixed and the formal fiber of $R$ at $\frak p$ is Cohen-Macaulay.
\end{proposition}

\begin{proof} Set $t=\dim (R/\frak p).$  From our assumption, we get by Lemma \ref{L:1a} and Lemma \ref{L:1}(ii),(iv)  that
$$\Rad (\frak p\widehat R)=\Rad\big((\Ann_RH^t_{\frak m}(R/\frak p))\widehat R\big)=\Rad \big(\Ann_{\widehat R}H^t_{\frak m}(R/\frak p)\big)=\Rad(U_{\widehat R}(\frak p\widehat R)).$$ So, $R/\frak p$ is  quasi unmixed. 

Assume that the formal fiber of $R$ at $\frak p$ is not Cohen-Macaulay. Then there exists prime ideal $\frak P$ of $\widehat R$ such that $\frak P\cap R=\frak p$ and $\widehat R_{\frak P}/\frak p \widehat R_{\frak P}$ is not Cohen-Macaulay. We choose such a $\frak P$ such that $\frak P$ is a minimal element of the non Cohen-Macaulay locus $\nCM (\widehat R/\frak p \widehat R)$ of $\widehat R/\frak p \widehat R$. Hence $\nCM (\widehat R_{\frak P}/\frak p\widehat R_{\frak P})=\{\frak P R_{\frak P}\}.$
Let $\frak Q \widehat R_{\frak P} \in \Ass_{\widehat R_{\frak P}}(\widehat R_{\frak P}/\frak p \widehat R_{\frak P})$, where  $\frak Q\in \Ass(\widehat R/\frak p \widehat R)$ and $\frak Q\subseteq \frak P$, $\frak Q\neq \frak P$. Then $\widehat R_{\frak Q}/\frak p \widehat R_{\frak Q}$ is Cohen-Macaulay and $\frak Q \widehat R_{\frak Q} \in \Ass(\widehat R_{\frak Q}/\frak p \widehat R_{\frak Q}).$ Hence,  $\ell(\widehat R_{\frak Q}/\frak p \widehat R_{\frak Q}) <\infty$ and hence $\frak Q\in \min \Ass(\widehat R/\frak p \widehat R).$ Since $R/\frak p$ is quasi unmixed as in the above proof, $\dim (\widehat R/\frak Q) =\dim (R/\frak p)$. We have
\begin{align*}
	\dim(\widehat R_{\frak P}/\frak Q \widehat R_{\frak P}) &=\dim(\widehat R/\frak Q)-\dim(\widehat R/\frak P)
	=\dim(R/\frak p)-\dim(\widehat R/\frak P)\\
	&= \dim (\widehat R/\frak p \widehat R)  -\dim(\widehat R/\frak P)
	=\dim (\widehat R_{\frak P}/\frak p \widehat R_{\frak P}).
\end{align*}
  Therefore, $\Assh(\widehat R_{\frak P}/\frak p \widehat R_{\frak P})=\Ass(\widehat R_{\frak P}/\frak p \widehat R_{\frak P})\setminus \{ \frak P \widehat R_{\frak P} \}.$ As  $\nCM (\widehat R_{\frak P}/\frak p\widehat R_{\frak P})=\{\frak P R_{\frak P}\}$, it follows that $\widehat R_{\frak P}/\frak p \widehat R_{\frak P}$ is generalized Cohen-Macaulay.
  Set $r=\dim (\widehat R_{\frak P}/\frak p \widehat R_{\frak P}).$ Since  $\widehat R_{\frak P}/\frak p \widehat R_{\frak P}$ is generalized Cohen-Macaulay and it is not Cohen-Macaulay, there exists an integer $i\leq r-1$ such that 
$0 < \ell_{\widehat R_{\frak P}} (H_{\frak P \widehat R_{\frak P}}^{i}(\widehat R_{\frak P}/\frak p \widehat R_{\frak P}))<\infty.$ Therefore, $\Att_{\widehat R_{\frak P}}(H_{\frak P \widehat R_{\frak P}}^{i}(\widehat R_{\frak P}/\frak p \widehat R_{\frak P}) )=\{\frak P \widehat R_{\frak P}\}$ by Lemma \ref{L:1}(i).  Note that $ \frak P \in \Att_{\widehat R} H_{\frak m}^{i+\dim(\widehat R/\frak P)}(\widehat R/\frak p \widehat R)$ by Weak general Shifted Localization Principle (see \cite[11.3.8]{BS}). Hence $\frak p\in \Att_R(H_{\frak m}^{i+\dim(\widehat R/\frak P)}(R/\frak p))$ by Lemma \ref{L:1}(iii). Therefore, we get by \cite[Corollary 11.3.5]{BS} and Lemma \ref{L:1}(ii) that
  \begin{align*}
  	 \dim \Big(\widehat R/\Rad\big (\Ann_{\widehat R}&H_{\frak m}^{i+\dim(\widehat R/\frak P)}(R/\frak p)\big)\Big) \leq i+\dim(\widehat R/\frak P)\\
  	 &<r +\dim(\widehat R/\frak P)
  	 =\dim(\widehat R/\frak p\widehat R)
  	 =\dim (R/\frak p) \\
  	 &\leq \dim \big(R/\Ann_RH_{\frak m}^{i+\dim(\widehat R/\frak P)}(R/\frak p)\big)\\
  	 &=\dim \Big(\widehat R/\Rad\big((\Ann_RH_{\frak m}^{i+\dim(\widehat R/\frak P)}(R/\frak p))\widehat R\big)\Big).
  \end{align*}
  So, $ \Rad \big({\mathrm{Ann}_{\widehat{R}} H^{i+\dim(\widehat R/\frak P)}_{\frak m}(R/\frak p)\big) \neq  \Rad\big(\Ann_{R} H^{i+\dim(\widehat R/\frak P)}_{\frak m}(R/\frak p)\widehat{R}}\big)$. This gives a contradiction to the assumption. Thus, the formal fiber of $R$ at $\frak p$ is Cohen-Macaulay.
  \end{proof}

For each integer $i\geq 0$ and each finitely generated $R$-module $M$ with $\dim_R(M)=d$, we set $\frak a_i(M)=\Ann_RH^i_{\frak m}(M)$ and $\frak a(M)=\frak a_0(M)\ldots  \frak a_{d-1}(M).$ Following N. T. Cuong \cite{C}, a system of parameters $x_1, \ldots , x_d$ of $M$ is called a {\it p-standard system of parameters} if $x_d\in\frak a(M)$ and $x_i\in\frak a(M/(x_{i+1}, \ldots , x_d)M)$ for all $i=1, \ldots , d-1.$ Following \cite[Theorems 1.2, 1.3]{CC},  $R$ is a quotient of a Cohen-Macaulay local ring if and only if every finitely generated $R$-module admits a  p-standard system of parameters, if and only if  $R$ admits a  p-standard system of parameters.

Now we are ready to prove Theorem \ref{T:3a}.

\begin{proof}[{\bf Proof of Theorem \ref{T:3a}}] (i) $\Rightarrow$ (ii). Let $M$ be a finitely generated $R$-module, $\frak p\in \Supp_R(M)$ and $i\geq 0$. We have by the assumption (i) and \cite[Theorem 1.1]{NQ} that   
	$$ \Att_{R_\frak p} (H^{i-\dim (R/\frak p)}_{\frak p R_{\frak p}}(M_{\frak p}))=\{\frak q R_{\frak p} \mid \frak q \in \Att_R H^i_{\frak m}(M), \frak q \subseteq \frak p\}. $$
 Then, by Lemma \ref{L:1}(ii) we get
	\begin{align*}\Rad\Big(\Ann_{R_\frak p} H^{i-\dim (R/\frak p)}_{\frak p R_{\frak p}}(M_{\frak p}&) \Big)=\underset {\frak q \in \Att_R H^i_{\frak m}(M), \frak q \subseteq \frak p}{\bigcap\frak q R_{\frak p}}=\Big(\underset{\frak q \in \Att_R H^i_{\frak m}(M),\, \frak q \subseteq \frak p}{\bigcap\frak q }\Big)R_{\frak p}\\
&=\Big(\underset{\frak q \in \Att_R H^i_{\frak m}(M)}{\bigcap\frak q }\Big)R_{\frak p}=\big(\Rad(\Ann_R H^i_{\frak m}(M))\big)R_{\frak p}\\
&=\Rad\Big((\Ann_R H^i_{\frak m}(M))R_{\frak p}\Big)\end{align*}

(ii) $\Rightarrow$ (i). Let $M$ be a finitely generated $R$-module.  For any integer $i\geq 0$, we claim that $\dim\big(R/\Ann_RH^i_{\frak m}(M)\big)\leq i$. Indeed, assume in contrary that $\dim\big(R/\Ann_RH^i_{\frak m}(M)\big)>i$ for some $i$. By Lemma \ref{L:1}(ii), there exists $\frak p\in\Att_R(H^i_{\frak m}(M))$ such that $\dim(R/\frak p)>i$. Hence $\frak p\supseteq \Ann_RH^i_{\frak m}(M)$ by Lemma \ref{L:1}(ii)  and $H^{i-\dim(R/\frak p)}_{\frak p R_{\frak p}}(M_{\frak p})=0.$ So, 
$$\Rad\Big(\big(\Ann_RH^i_{\frak m}(M)\big)R_{\frak p}\Big)\neq R_{\frak p}=\Rad\Big(\Ann_{R_{\frak p}}H^{i-\dim(R/\frak p)}_{\frak p R_{\frak p}}(M_{\frak p})\Big).$$  This gives a contradiction to the assumption (ii). So the claim is proved. From the definition of p-standard system of parameters, it follows by the claim that every finitely generated $R$-module admits a p-standard system of parameters. Therefore, $R$ is a quotient of a Cohen-Macaulay local ring.

(i)$\Rightarrow$(iii). We get by assumption (i) and \cite[Theorem 1.1]{NQ} that 
	$$\mathrm{Att}_{\widehat{R}}H^i_{\frak m}(M)=\bigcup_{\frak p \in \mathrm{Att}_RH^i_{\frak m}(M)} \mathrm{Ass}_{\widehat{R}}(\widehat{R}/\frak p\widehat{R}).$$
	Therefore, by Lemma \ref{L:1}(ii) we have
	\begin{align*}\Rad\Big({\mathrm{Ann}_{\widehat{R}} H^i_{\frak m}(M)}\Big) &= \bigcap_{\frak p \in \mathrm{Att}_RH^i_{\frak m}(M) } \Big(\bigcap_{\frak P \in \mathrm{Ass}_{\widehat{R}}\widehat{R}/\frak p\widehat{R}}\frak P\Big)\\
	&= \bigcap_{\frak p \in  \mathrm{Att}_RH^i_{\frak m}(M) } \Rad(\frak p\widehat{R})
		= \Rad\Big(\big(\bigcap_{\frak p \in \mathrm{Att}_RH^i_{\frak m}(M)}\frak p\big)\widehat{R}\Big)\\
		&= \Rad\Big(\big(\Rad(\mathrm{Ann}_{R} H^i_{\frak m}(M))\big)\widehat{R}\Big)\\
		&= \Rad\big(\big(\mathrm{Ann}_{R} H^i_{\frak m}(M)\big)\widehat{R}\big).
	\end{align*}

(iii)$\Rightarrow$(i).  Let $\frak p\in\Spec (R).$ We get by assumption (iii) and by Proposition \ref{P:1c} that $R/\frak p$ is quasi unmixed and the formal fiber of $R$ at $\frak p$ is Cohen-Macaulay. Therefore, $R$ is universally catenary and all of its formal fibers are Cohen-Macaulay. By \cite[Corollary 1.2]{Kaw}, $R$ is a quotient of a Cohen-Macaulay local ring. 
\end{proof}

 Note that if $R$ is a quotient of a Cohen-Macaulay local ring then $R/\frak p$ is unmixed for all $\frak p\in\Spec(R)$. The converse statement is not true. M. Brodmann constructed a Noetherian local domain $(R,\frak m)$ such that $R/\frak p$ is unmixed for all $\frak p\in\Spec(R)$, but $R$ is not a quotient of a Cohen-Macaulay local ring. In this case, the annihilator of all top local cohomology modules $H^{\dim_R(M)}_{\frak m}(M)$ is compatible under localization and completion by Theorems \ref{T:1a}, \ref{T:2a}, but there exists a finitely generated $R$-module $M$ and an integer $i<\dim_R(M)$ such that the annihilator  of $H^i_{\frak m}(M)$ is neither compatible under localization nor compatible under completion.  

Now we give some examples to clarify the results of this paper.

\begin{example} \label{E:1} {\rm Let $(R,\frak m)$ be the Noetherian local domain of dimension $2$ constructed by D. Ferrand and M. Raynaud \cite{FR} such that $\dim(\widehat R/\frak P)=1$ for some $\frak P\in\Ass (\widehat R).$  We examine the compatibility of the annihilator of $H^{2}_{\frak m}(R)$ and $H^{1}_{\frak m}(R)$ under localization and completion:
		
(i) Since $R$ is catenary, the annihilator of $H^2_{\frak m}(R)$ is compatible under localization by Theorem \ref{T:1a}. Since $R$ is not unmixed, $U_{\widehat R}(0) \neq 0$. It follows from Lemma \ref{L:1a} that $\Ann_{\widehat R} H_{\frak m}^{2}(R)=U_{\widehat R}(0)$. As $R$ is a domain, $U_R(0)=0$, hence $\big(\Ann_R H_{\frak m}^{2}(R)\big)\widehat R=0$ by Lemma \ref{L:1a}. Thus, the annihilator of $H^{2}_{\frak m}(R)$ is not compatible under completion.

(ii) As $\frak P \in \Ass (\widehat{R})$ and $R$ is a domain, $\frak P \cap R=0$. For $\frak p=0$, we have $H^{1-\dim (R/\frak p)}_{\frak p R_{\frak p}}(R_{\frak p})=0$, hence  $\Ann_{R_{\frak p}}H^{i-\dim (R/\frak p)}_{\frak p R_{\frak p}}(R_{\frak p})=R_{\frak p}.$ It follows by \cite[11.3.3]{BS} that $\frak P \in \Att_{\widehat R} H_{\frak m}^1(R)$. So, $\frak P\cap R =0\in \Att_R H_{\frak m}^1(R)$ and hence $\Ann_{R}H^{1}_{\frak m}(R)=0$ by Lemma \ref{L:1}(ii),(iii). Thus,  $$(\Ann_{R}H^{1}_{\frak m}(R))_{\frak p}=0\neq R_{\frak p}=\Ann_{R_{\frak p}}H^{i-\dim (R/\frak p)}_{\frak p R_{\frak p}}(R_{\frak p}),$$ i.e. the annihilator of $H^{1}_{\frak m}(R)$ is not compatible under localization.  It follows by Lemma \ref{L:1}(ii) and \cite[11.3.5]{BS} that $\dim \big(\widehat R/ \Ann_{\widehat R}H^1_{\frak m}(R)\big)=1$. Hence
$$\Ann_{\widehat R}H^1_{\frak m}(R)\neq 0=\big(\Ann_{R}H^{1}_{\frak m}(R)\big)\widehat R,$$ i.e.  the annihilator of $H^1_{\frak m}(R)$ is not compatible under completion.
}\end{example} 

For local domains of higher dimension, we consider the following example.

\begin{example} \label{E:2} {\rm  Let $K$ be a field. Let $K[[x_1, x_2, \ldots, x_{d+1}]]$ be the ring of formal power series in $d+1$ variables with coefficients in $K$. By \cite{CL}, there exists a Noetherian local domain $(R,\frak m)$  such that  $\widehat R=K[[x_1, x_2, \ldots, x_{d+1}]]/(x_1^2)\cap (x_1^3,x_2),$ see also \cite[Example 3.8]{CNN}.  Then $\dim (R)=d$ and $R$ is not unmixed, but $R$ is quasi unmixed and hence $R$ is catenary. We examine the compatibility of the annihilator of $H^{d}_{\frak m}(R)$ and $H^{d-1}_{\frak m}(R)$ under localization and completion:

(i) The annihilator of $H^{d}_{\frak m}(R)$ is compatible under localization by Theorem \ref{T:1a}. The annihilator of $H^{d}_{\frak m}(R)$ is not compatible under completion since  $\Ann_{\widehat{R}} H^d_{\frak m}(R)=(x_1^2)$ and  $\big(\Ann_{R} H^d_{\frak m}(R)\big)\widehat{R}=(x_1^2)\cap (x_1^3, x_2)$ by Lemma \ref{L:1a}. Note that 
$$\Rad \big(\Ann_{\widehat{R}} H^d_{\frak m}(R)\big) = \Rad\big(\big(\Ann_{R} H^d_{\frak m}(R)\big)\widehat{R}\big)=(x_1).$$

(ii) The annihilator of $H^{d-1}_{\frak m}(R)$ is not compatible under completion. It is clear that $\frak P:=(x_1, x_2)\in \Ass(\widehat R)$. Hence  $\frak P\in \Att_{\widehat R} H_{\frak m}^{d-1}(R)$ by \cite[11.3.3]{BS}. So, we get by Lemma \ref{L:1}(ii) and \cite[11.3.5]{BS} that $\dim \big(\widehat R/\Ann_{\widehat R} H_{\frak m}^{d-1}(R)\big)=d-1$. Since $R$ is a domain, $0=\frak P\cap R\in  \Att_R H_{\frak m}^{d-1}(R)$ by Lemma \ref{L:1}(iii). So, $\Ann_RH_{\frak m}^{d-1}(R)=0$ and hence 
$$\Ann_{\widehat R} H^{d-1}_{\frak m}(R)\neq 0=\big(\Ann_R H^{d-1}_{\frak m}(R)\big)\widehat{R}.$$
The annihilator of $H^{d-1}_{\frak m}(R)$ is not compatible under localization. With $\frak p=0$ we have $\dim(R/\frak p)=d$ and
$$\Ann_{R_{\frak p}} H^{(d-1)-\dim (R/\frak p)}_{\frak p R_{\frak p}}(R_{\frak p})=R_{\frak p}\neq 0=\big(\Ann_R H^{d-1}_{\frak m}(R)\big)R_{\frak p}.$$}
\end{example}

\end{document}